\documentclass[review]{elsarticle}

\usepackage{lineno,hyperref}
\usepackage{xcolor}
\modulolinenumbers[5]
\usepackage{amssymb}
\usepackage{amsmath}
\usepackage{amsthm}
  \usepackage{lineno}

\newcommand{\R}{\mathbb{R}}

\newcommand{\epi}{\mathrm{epigraph}}
\newcommand{\hyp}{\mathrm{hypograph}}  %epigraph
\newcommand{\grap}{\mathrm{graph}}
  \newtheorem{thm}{Theorem}%[section]
\newtheorem{lem}{Lemma}%[section]
\newtheorem{cor}{Corollary}%[section]
\newtheorem{exam}{Example}%[section]
\newtheorem{defi}{Definition}%[section]
\begin{document}

\begin{frontmatter}

\title{A new separable property of the joint numerical range of quadratic functions and its applications to the 
Smallest Enclosing Ball Problem }

\author[mymainaddress]{Van-Bong Nguyen\corref{mycorrespondingauthor}}
\cortext[mycorrespondingauthor]{Corresponding author}
\ead{nvbong@ttn.edu.vn}
\address[mymainaddress]{Department of Mathematics, Tay Nguyen
University, 632090, Vietnam}
\author[mysecondaddress]{Huu-Quang Nguyen}
\ead{quangdhv@gmail.com}

\address[mysecondaddress]{Department of Mathematics, Vinh
University, 463360, Vietnam}

\begin{abstract}
We explore separable property of the joint numerical range $G(\R^n)$ of a special class of quadratic functions and apply it to solving the smallest 
enclosing ball (SEB) problem which asks to find a ball $B(a,r)$ in  $\Bbb R^n$ with smallest radius $r$
such that $B(a,r)$ contains the intersection $\cap_{i=1}^mB(a_i,r_i)$ of $m$ given balls $B(a_i,r_i).$
We show that $G(\R^n)$ is convex if and only if ${\rm rank}\{a_1-a, a_2-a, \ldots, a_m-a\}\le n-1.$  Otherwise,
${\rm rank}\{a_1-a, a_2-a, \ldots, a_m-a\}=n$ and $G(\R^n)$ is not convex. In this case we propose a new set 
$G(\R^n)^\bullet$ which allows to show that if $m=n$ then $G(\R^n)^\bullet$ is convex even $G(\R^n)$ is not. Importantly,
the separable property of $G(\R^n)^\bullet$ then implies the separable property for $G(\R^n).$ As a result,
a new progress on solving the SEB problem  is obtained. 

\end{abstract}


\end{frontmatter}

%\linenumbers

\section{Introduction}\label{sec:intro}
Let $g(x)=x^Tx-2a^Tx+\theta, g_i(x)=x^Tx-2a_i^Tx+\theta_i, i=1,2,\ldots,m,$ 
where $a, a_i\in\Bbb R^n; \theta, \theta_i\in\Bbb R,$ be strictly convex quadratic functions and denote
$$\Lambda=\{z=(z_0, z_1, z_2, \ldots, z_{m})^T: z_0<0, z_i\leq 0, i=1,2,\ldots, m\}\subset \R^{m+1}$$  the negative orthant in
$\Bbb R^{m+1}.$  
In this paper we study the separable property of the 
joint numerical range 
$$G(\R^n)=\{(-g(x), g_1(x), \ldots, g_m(x))^T, x\in\R^n\}\subset \R^{m+1}$$  of the quadratic mapping
$G=(-g,g_1, \ldots, g_m):\R^n\longrightarrow \R^{m+1}$
 in the sense that we find conditions under which  $G(\R^n)$ and $\Lambda$ can be separated by a hyperplane whenever
$G(\R^n)\cap \Lambda=\emptyset.$ Our study has the main motivation from the smallest enclosing ball (SEB) problem which asks to find
a ball $B(a, r)$ centered at $a\in\Bbb R^n$ with smallest radius $r$   such that $B(a,r)$ contains the intersection 
$\cap_{i=1}^mB(a_i,r_i)$
of given balls $B(a_i, r_i), i=1,2,\ldots,m.$ To see the relation between the separable property of  $G(\R^n)$ and the SEB problem,
we make a short review as follows. 
 The balls can be expressed as:
$B(a_i,r_i)=\{x\in\Bbb R^n: g_i(x)\le0 \}, i=1,2,\ldots,m,$ and $B(a, r)=\{x\in\Bbb R^n: g(x)\le0\}.$ Suppose that the intersection $\cap_{i=1}^mB(a_i,r_i)$
has a nonempty interior.  The SEB problem is then
restated in the following optimization problem  
\begin{align}\label{problem}
	\begin{array}{llll}
		&\min & r& \\
		&{\rm s.t.}& \cap_{i=1}^mB(a_i,r_i)\subset B(a,r),\\
		& &r>0, a\in\R^n,
	\end{array}
\end{align}
 where $r=\sqrt{\|a\|^2-\theta}$ is  needed to find and $r_i=\sqrt{\|a_i\|^2-\theta_i}$ are given.
  The constraint $\cap_{i=1}^mB(a_i,r_i)\subset B(a,r)$  means that 
  $(x\in\Bbb R^n) ~g_1(x)\le 0, \ldots, g_m(x)\le 0 \Rightarrow g(x)\le0.$ Or equivalently,
 the system of inequalities $ g_1(x)\le 0, \ldots, g_m(x)\le 0 , - g(x)<0$ has no root, and thus  $G(\R^n)\cap\Lambda=\emptyset.$ 
 Now, if $G(\R^n)$ and $\Lambda$ are separable by a hyperplane, 
problem \eqref{problem} is reformable as an SDP problem and thus solvable in polynomial time \cite{Polik-Terlaky07}. So the question is that under what conditions, 
$G(\R^n)$ and $\Lambda$ are separable by a hyperplane?   

Separable property of the joint numerical range is an interesting topic in literature with various applications. Normally, if 
the joint numerical range is convex, it is separable. Dines \cite{Dines41} showed
 that  the joint numerical range of two homogeneous quadratic functions  $W(A,B)=\{(x^TAx, x^TBx), x\in\Bbb R^n\}\subset\Bbb R^2$ is always convex, where $A, B\in\Bbb R^{n\times n}$ are symmetric matrices. 
 This result was then applied by
 Yakubovich to obtain the well known S-Lemma \cite{Yakubovich}.
   Brickman \cite{Brickman} showed that the restriction of $W(A,B)$ on the unit sphere: 
 $\{(x^TAx, x^TBx),\|x\|=1\}\subset\Bbb R^2$ with $n\ge3$ is still convex. This helps to obtain an extension of the S-Lemma, please see \cite{Polik-Terlaky07,Bong}. However, for nonhomogeneous quadratic functions or more than two homogeneous quadratic functions, the convexity
 may not be guaranteed. 
 Polyak \cite{Polyak98} showed that
 the joint numerical range of three quadratic forms $W(A,B,C)=\{(x^TAx, x^TBx, x^TCx), x\in\Bbb R^n\}$
 with $n\ge3$  is convex and the quadratic forms $x^TAx, x^TBx, x^TCx$ has no common zero except zero  if and only if there is a positive definite combination of $A, B, C,$ i.e., there are scalars $\lambda, \gamma, \eta$ such
 that $\lambda A+\gamma B+\eta C\succ0.$
  For quadratic functions $q_i(x)=x^TA_ix+a_i^Tx+\alpha_i, i=1,2,$   
 the joint numerical range $W(q_1,q_2)=\{(q_1(x), q_2(x)), x\in\Bbb R^n\}$ with $n\ge2$ is convex if there are some
 $\lambda, \gamma\in\Bbb R$ such that $\lambda A_1+\gamma A_2\succ0$  \cite{Polyak98}.    
  For more results on convexity and separable property of the joint numerical range,
 the readers are encouraged to see the well known survey \cite{Polik-Terlaky07} and more recent studies
  \cite{FloresBazan-Opazo16,Quang-Sheu18,Quang-Chu-Sheu20,Flores-Bazan}  with
 references therein. 
 For the joint numerical range  $G(\R^n),$ a special type of quadratic mapping,  Beck \cite{Beck07} showed that  if $m\le n-1$ then $G(\R^n)$ is convex.

The model \eqref{problem}  is  a special case of the outer approximation problems
which find the minimum volume ellipsoid enclosing a given set $S$ in $\R^n,$ see for example
 \cite{Boyd-Ghaoui-Feron-Balakrishnan94,Vandenberghe} for the models of the problem and \cite{Berg,Eliosoff} for some applications  such as the applications in  environmental science, in pattern
recognition, in protein analysis, in analysis of party spectra, of
 the problem.
 Depending on the set $S,$ the outer approximation problem  is known tractable or not. For example, if $S$ is a finite set of points, a polyhedron, or
  a union of ellipsoids, the problem is solvable in polynomial time
\cite{Boyd-Ghaoui-Feron-Balakrishnan94,Vandenberghe0,Beck07}. Algorithms for those cases can be found, for example in  \cite{Xu} for the case
$S$ is a union of circles  in the Euclidian plane;  in
  \cite{Zhou} when $S$ is a union of balls in $\R^n;$ in
 \cite{Yildirim} when $S$ is a finite set of
points in $\R^n.$
However, when $S$ is the intersection of ellipsoids, the problem is considered difficult, because
only checking whether an ellipsoid contains the intersection of given ellipsoids is already NP-hard
\cite[Sect. 3.7.2]{Boyd-Ghaoui-Feron-Balakrishnan94}. For problem \eqref{problem},
 $S$ is the intersection of balls and the problem is  also known as the ball version of the
outer approximation problems. Beck \cite{Beck07, Beck09} showed that if $m\le n$ then the SEB problem \eqref{problem} is solvable
in polynomial time.  On the other hand, a
recent study by Xia et al. \cite{Xia-Yang-Wang20}  showed that the problem is in general NP-hard.
Those results raise a question that how large the  number $m$ of given balls can be such that the SEB problem is still solvable in polynomial time?

In this paper  
we show that $G(\R^n)$ is convex if and only if ${\rm rank}\{a_1, a_2, \ldots, a_m\}\le n-1.$    
%This slightly extends 
%Beck's result \cite[Theorem 3.2]{Beck07} in the sense that
%the number $m$ of given balls may be any provided that the rank of the vectors $a_1, a_2, \ldots, a_m$ is not greater than $n-1.$
   Otherwise, ${\rm rank}\{a_1, a_2, \ldots, a_m\}=n$ and  $G(\R^n)$ is then not convex. 
   In this case, we propose a new set denoted by
    $G(\R^n)^\bullet$  to reveal several interesting properties of
$G(\R^n),$ please see Definition \ref{de} below. We show that if
${\rm rank}\{a_1, a_2, \ldots, a_m\}=n$ and $m=n,$  then $G(\R^n)^\bullet$ is convex even $G(\R^n)$ is not;  especially,
$G(\R^n)\cap\Lambda=\emptyset$ implies $G(\R^n)^\bullet\cap\Lambda=\emptyset;$ we note that in this case $G(\R^n)\subsetneq G(\R^n)^\bullet.$
This new separable property of $G(\R^n)$ allows to obtain a new extension of the S-Lemma, please see Corollary  \ref{thm2} below, which helps 
     to reformulate the SEB
problem \eqref{problem} as an SDP problem. Then, using the same analysis technique in \cite{Beck07} we further reformulate the resulting SDP problem as a
 problem of minimizing a convex quadratic function over the unit simplex.
 %Our results together with Xia's result in \cite{Xia-Yang-Wang20} naturally raise a conjecture of whether the SEB problem is NP-hard if and only if
%${\rm rank}\{a_1, a_2, \ldots, a_m\}=n$ and $m>n?$
\section{On the convexity of $G(\R^n)$}
 In this section we propose a necessary and sufficient condition for $G(\R^n)$ to be convex. Our proof is based on 
  the following result.

\begin{lem}[\cite{Beck07}]\label{L1}
If  $m\le n-1,$  the set
$$S=\{(x^Tx-2a^Tx, a_1^Tx, a_2^Tx, \ldots, a_m^Tx)^T, x\in\Bbb R^n\}\subset\Bbb R^{m+1}$$
is closed and convex.
\end{lem}
%Using Lemma \ref{L1} we prove the main result of this subsection as follows.

%  
 The main result of this section is now obtained as follows.
\begin{thm}\label{thm1}
$G(\R^n)$   is convex if and only if ${\rm rank}\{a_1-a, a_2-a, \cdots, a_m-a\}\leq n-1.$   
 \end{thm}
\proof
Suppose first that ${\rm rank}\{a_1-a, a_2-a, \cdots, a_m-a\}\le n-1,$  we shall show that $G(\R^n)$  is convex. Let  
$p={\rm rank}\{a_1-a, a_2-a, \cdots, a_m-a\}$ and suppose without loss of generality 
that the first $p$ vectors are linear independent, i.e., ${\rm rank}\{a_1-a, a_2-a, \ldots, a_p-a\}=p\le n-1.$ 
 Then  for each $1\le i\le m, a_i-a$ is uniquely written  as a linear combination
of $a_1-a, \cdots, a_p-a$ as follows  
\begin{equation}\label{234qqq}
	a_i-a=\sum_{j=1}^{p}{\alpha_{ij}}(a_j-a), ~ i=1,2, \ldots, m,
\end{equation}
where $\alpha_{ij}\in\Bbb R.$ Let $C\subset \R^{p+1}$ be a set defined by 
$$C:=\left\{\left(-x^Tx+2a^Tx, -2(a_1-a)^Tx, \ldots, -2(a_p-a)^Tx \right)^T, x\in\Bbb R^n\right\}\subset \R^{p+1}.$$
Using Lemma \ref{L1}, we can easily show that $C$
is closed and convex. Let $T: \R^{p+1} \rightarrow \R^{m+1}$ be a linear mapping
defined by 
$$T(z_0, z_1 \ldots, z_p) = \left(z_0, \sum_{j=1}^p \alpha_{1j}z_j, \sum_{j=1}^p \alpha_{2j}z_j,\ldots, \sum_{j=1}^p \alpha_{mj}z_j\right)^T.$$
We then have
$$
	 T(C) = \left\{\left(-x^Tx+2a^Tx, -2(a_1-a)^Tx, \ldots, -2(a_p-a)^Tx,  \ldots, -2(a_m-a)^Tx\right)^T\right\}.
$$
 Let    $L: \R^{m+1} \rightarrow \R^{m+1}$  be an affine transformation defined by  
\begin{align*} 
& L(z_0, z_1 \ldots, z_m)=(z_0, z_1-z_0, \ldots, z_m-z_0)^T+(-\theta, \theta_1,\ldots,\theta_m)^T.
\end{align*}
Then we observe that
\begin{align*} 
  L(T(C))=&\left\{ \left(-x^Tx+2a^Tx-\theta, x^Tx-2a_1^Tx+\theta_1,   \ldots, x^Tx-2a_m^Tx+ \theta_m\right)^T, x\in\Bbb R^n\right\} \\
 =&G(\Bbb R^n).
\end{align*}

 Since 
  the union $L\circ T$ is an affine mapping and
    $C$ is   convex, its image
$(L\circ T)(C)$ is convex,  $G(\Bbb R^n)$ is  thus convex.

Conversely, we suppose on contrary that  $G(\Bbb R^n)$ is convex but  ${\rm rank}\{a_1-a, a_2-a, \ldots, a_m-a\}=n.$
 Let $D$ be defined by
$$D=\left\{\left(-x^Tx+2a^Tx, -2(a_1-a)^Tx, \ldots, -2(a_m-a)^Tx \right)^T, x\in\Bbb R^n\right\}\subset \R^{m+1}.$$ 
Applying Theorem 2.2 in  \cite{Beck07}, we can show that $D$ is not convex. But now, by defining the following affine transformation
 $F:\R^{m+1}\longrightarrow\R^{m+1},$  
 $$F(z_0, z_1, \ldots, z_m)=(z_0, z_1+z_0, \ldots, z_m+z_0)+(\theta, -(\theta_1-\theta), \ldots, -(\theta_m-\theta))$$
we will have
 $$F(G(\R^n))=\{(-x^Tx+2a^Tx, -2(a_1-a)^Tx, \ldots, -2(a_m-a)^Tx)\}=D.$$
That is, $D$ is the image via an affine transformation of the convex set $G(\R^n),$ it must be convex and we get a contradiction.
So it must be that
${\rm rank}\{a_1-a, a_2-a, \ldots, a_m-a\}<n.$
 \endproof
%\begin{cor}\label{cor1}
%If ${\rm rank}\{a_1-a, a_2-a, \ldots, a_m-a\}=n$  then $G(\R^n)$ is not convex.  
%\end{cor}

%The following example illustrates Corollary \ref{cor1}.
 
\begin{exam}\rm
	Let $g(x)= x_1^2+x_2^2-2x_1-2x_2,$ $g_1(x)=x_1^2+x_2^2-2x_2,$ $g_2(x)=x_1^2+x_2^2-2x_1$
be three convex quadratic functions in $\Bbb R^2.$
	Here, $a=(1,1)^T, a_1=(0, 1)^T, a_2=(1,0)^T$ and ${\rm rank}\{a_1-a, a_2-a\}=2=n$
and $\theta=\theta_1=\theta_2=0.$
	We check that 
$$G(\R^2)=\{(-g(x), g_1(x), g_2(x))^T\}\subset \Bbb R^3$$ 
is not convex. Indeed, let
	${u}=(1,0)^T, {v}=(0,1)^T\in\Bbb R^2$ we have $G(u)=(1,  1,-1)^T, G(v)=(1, -1,1)^T\in G(\R^2).$ The
	midpoint of the line segment $[G(u), G(v)]$ is $(1,0,0)^T.$ However, there is no point $x\in\Bbb R^2$ such that
$(-g(x), g_1(x), g_2(x))^T=(1,0,0)^T.$ That is the line segment $[G(u), G(v)]$ is not contained in $G(\R^n)$ and  so $G(\R^n)$ is not convex.
\end{exam}

%\subsection{When $G(\R^n)$ is not convex   }\label{sect3}
\section{Separable property of $G(\R^n)$}

%\section{A new separable property of the image $G(\R^n)$ }
As mentioned, if $G(\R^n)$ is convex, the separable property holds for it. By Theorem \ref{thm1},
$G(\R^n)$ is convex if and only if  ${\rm rank}\{a_1-a, a_2-a, \ldots, a_m-a\}\le n-1.$ In this section, 
we will show that when ${\rm rank}\{a_1-a, a_2-a, \ldots, a_m-a\}=n=m,$
the separable property still holds for $G(\R^n)$ even it is not convex. To do so, we need first to give the following new concept.
  
\begin{defi}\label{de}
	For a set $D\subset \R^n,$ we define a set $D^{\bullet}\subset\Bbb R^n$ as follows: 
$$D^{\bullet}=\{x\in \Bbb R^n: x=\lambda u+(1-\lambda)v; ~   u, v\in D,   0\le\lambda \le 1 \}. $$  
\end{defi}
By this definition, $D^{\bullet}$  is the set of the convex combinations of points of $D$ and  looks like the convex hull ${\rm conv}(D)$ of $D.$ But, in fact, 
$D^{\bullet}$ is different from  ${\rm conv}(D).$   
 Indeed, we consider the following
  example: let $D=\{u=(0,0)^T, v=(1,0)^T, \omega=(0,1)^T\}\subset \R^2$ be the set of three points $u, v, \omega$ in $\Bbb R^2.$
Then $D^{\bullet}$ is  the three edges of the triangle formed by the three vertices $u, v, \omega,$   
while ${\rm conv}(D)$ is the triangle  of those three vertices. More concrete,  
for any $D\subset\Bbb R^n,$  it holds that
	$$D\subset D^{\bullet}\subset{\rm conv}(D).$$
The following two results are important but their proofs are so simple, we thus omit them.
\begin{lem}\label{LL}
 $D=D^{\bullet}={\rm conv}(D)$
if  and only if $D$ is convex.  
\end{lem}
	
\begin{lem}\label{B0}
	$D^{\bullet}$ is invariant under a nonsingular affine transformation. That is, if
	$L$ is an invertible affine transformation in $\Bbb R^n$ then
	$$L(D^{\bullet})=L(D)^\bullet.$$
\end{lem}

 %
%	
% 
% The following property is important to our later analysis.
\begin{lem}\label{thm21qq}	
\begin{itemize}
\item[(i)]	If $f(x)=x^TAx-2a^Tx+\theta$ is a  convex quadratic function, i.e.,  $A\succeq0$ and $A\ne0,$ then $ \grap (f)^{\bullet}=\epi(f);$
\item[(ii)]	If $f(x)=x^TAx-2a^Tx+\theta$ is a  concave quadratic function, i.e.,  $A\preceq0$ and $A\ne0,$ then $ \grap (f)^{\bullet}=\hyp(f).$
\end{itemize}
\end{lem}
\proof
We prove only (i). The proof of (ii) is done similarly.

We first show  $\grap (f)^{\bullet}\subset\epi(f).$ Let $\omega\in \grap (f)^{\bullet},$
by definition, there will be $u=(x_1, y_1), v=(x_2,y_2)\in \grap (f), y_1=f(x_1), y_2=f(x_2),$ such that
$\omega=\lambda u+(1-\lambda)v=(\lambda x_1+(1-\lambda)x_2, \lambda y_1+(1-\lambda)y_2)$
for some $0\le \lambda\le 1.$ Since $f$ is convex, we have
$$f(\lambda x_1+(1-\lambda)x_2)\le \lambda f(x_1)+(1-\lambda)f(x_2)=\lambda y_1+(1-\lambda)y_2.$$
This shows   $\omega\in \epi(f)$ and implies $\grap (f)^{\bullet}\subset\epi(f).$

For the converse,  let $\omega=(\bar{x},t)\in\epi(f)$ be any point in the epigraph of $f,$ i.e.,  $ f(\bar{x})\le t,$ we need to 
 show that $\omega\in \grap (f)^{\bullet}.$  Indeed,
 \begin{enumerate}
   \item If  $ f(\bar{x})= t$ then $\omega\in\grap(f)$ so  $\omega\in \grap (f)^{\bullet}.$ 
   \item  If  $ f(\bar{x})< t,$ we now show that there exist $u, v\in \grap (f) $ such that
$\omega=\lambda u+(1-\lambda)v,$  $0\le \lambda\le 1.$ Specifically, the points $u, v$ are of the form
$u=(\hat{x},t), v=(\tilde{x},t)$ such that $f(\hat{x})=f(\tilde{x})=t$ and $\bar{x}=\lambda\hat{x}+(1-\lambda)\tilde{x}.$
Observe that $\hat{x}$ and $\tilde{x}$ are solutions of the equation $f(x)=t.$ 
For simplicity when solving this equation we make $A$ diagonal as follows.  
Since $A\succeq0$ and $A\ne0,$  there is a nonsingular matrix $P$ such that $P^TAP={\rm diag}(\alpha_1, \alpha_2, \ldots, \alpha_n),$
  $\alpha_i\ge0,~\forall i=1,2,\ldots,n,$ and $\alpha_j>0$ for some 
$j.$   Let $x=Py,$ then 
$$x^TAx-2a^Tx+\theta=y^T(P^TAP)y-2(P^Ta)^Ty+\theta=\sum_{i=1}^n\alpha_iy_i^2-2\sum_{i=1}^na_iy_i+\theta,$$
where $y^T=(y_1, y_2, \ldots, y_n)$ and $(P^Ta)^T=(a_1,a_2, \ldots, a_n).$
 The equation $f(x)=t$ then becomes
 \begin{align}\label{eq1}
 \sum_{i=1}^n\alpha_iy_i^2-2\sum_{i=1}^na_iy_i+\theta=t.
 \end{align} 
 Let $\bar{y}=P^{-1}\bar{x}=(\bar{y}_1, \ldots, \bar{y}_{j-1}, \bar{y}_j, \bar{y}_{j+1},\ldots,\bar{y}_n)^T.$
 In \eqref{eq1} we fix   $y_i=\bar{y}_i$ for all $i\ne j,$ then \eqref{eq1}  becomes the equation of only one unknown $y_j.$
 Let 
 $$g(y_j)=\alpha_jy_j^2-2a_jy_j+\sum_{i\ne j} \alpha_i\bar{y}_i^2-2\sum_{i\ne j} a_i\bar{y}_i+\theta.$$
 Since $\alpha_j>0,$ $g(y_j)$ is a strictly convex function of $y_j,$ and  $g(\bar{y}_j)<t.$ Those imply that the equation 
 $g(y_j)=t$ has two distinct solutions denoted by $\hat{y}_j$ and $\tilde{y}_j$ such that
 $\bar{y}_j$ is in the line segment connecting two endpoints  $\hat{y}_j$ and $\tilde{y}_j.$ That is, there is a value
 $0\le \lambda\le 1$ such that $\bar{y}_j=\lambda\hat{y}_j+(1-\lambda)\tilde{y}_j.$
  Let  
  $$\hat{y}=(\bar{y}_1, \ldots, \bar{y}_{j-1}, \hat{y}_j, \bar{y}_{j+1},\ldots,\bar{y}_n)^T, \tilde{y}=(\bar{y}_1, \ldots, \bar{y}_{j-1}, \tilde{y}_j, \bar{y}_{j+1},\ldots,\bar{y}_n)^T.$$  
  We have $\bar{y} =\lambda\hat{y} +(1-\lambda)\tilde{y}$ and 
  $$\bar{x}=P\bar{y}=\lambda P\hat{y}+(1-\lambda)P\tilde{y}
  = \lambda\hat{x}+(1-\lambda)\tilde{x},$$
  where $\hat{x}=P\hat{y}$ and $\tilde{x}=P\tilde{y}.$ We now have
    $f(\hat{x})=f(\tilde{x})=t.$ Let $u=(\hat{x},t), v=(\tilde{x},t)$ then
  $u,v\in\grap(f)$ such that $\omega=\lambda u+(1-\lambda)v$ with $0\le \lambda\le 1$ as desired. 
So $\omega\in \grap (f)^{\bullet}$ and that $\epi(f)\subset \grap (f)^{\bullet}.$
 \end{enumerate}
\endproof

The main result of this section is now stated as follows.

\begin{thm}\label{lm222}	  
If either ${\rm rank}\{a_1-a, a_2-a, \ldots, a_m-a\}<n$ or
		${\rm rank}\{a_1-a, a_2-a, \ldots, a_m-a \}=n=m,$ then the following hold.
	\begin{enumerate}
	 	\item[(i)]   $G(\R^n)^\bullet$ is convex;
		\item[(ii)]    $G(\R^n)\cap \Lambda=\emptyset$ implies
		  $G(\R^n)^\bullet\cap \Lambda=\emptyset.$
	\end{enumerate}
\end{thm}

\proof
(i) If  ${\rm rank}\{a_1-a, a_2-a, \cdots, a_m-a\}<n$ then, by Theorem \ref{thm1}, $G(\R^n)$ is convex so, by Lemma \ref{LL}, $G(\R^n)^\bullet=G(\R^n)$ and thus
  $G(\R^n)^\bullet$ is  convex.
 If  ${\rm rank}\{a_1-a, a_2-a, \cdots, a_m-a\}=n$ and
$m=n,$
we define a linear transformation   $H:\R^{m+1}\rightarrow \Bbb R^{m+1}$  by
\begin{equation}\label{216pp1}
	H(z_0, z_1, \ldots, z_m)= (z_1+z_0, \ldots, z_m+z_0, -z_0)^T.
\end{equation}
Then
\begin{eqnarray}\label{217nnn1}
	H(G(x))&=&(g_1(x)-g(x), \ldots, g_m(x)-g(x), g(x))^T \label{217pp1}\\
	&=&(-2(a_1-a)^Tx+\theta_1-\theta, \ldots, -2(a_m-a)^Tx+\theta_m-\theta, x^Tx-2a^Tx+\theta)^T. \nonumber
\end{eqnarray}  %Hence
 
Let \begin{eqnarray}\label{236qqq}
	\begin{split}
		&y_i=-2(a_i-a)^Tx+\theta_i-\theta, i=1, \ldots, m,
	\end{split}
\end{eqnarray}
and $A$ be a square matrix defined by
    $A=\begin{bmatrix}-2(a_1-a)^T\\ \vdots \\ -2(a_m-a)^T\end{bmatrix}.$
Then  we have
$$y=(-2(a_1-a)^Tx+\theta_1-\theta, \ldots, -2(a_m-a)^Tx+\theta_m-\theta)^T=Ax+\hat{\theta},$$
  where $\hat{\theta}=(\theta_1-\theta, \theta_2-\theta,\ldots, \theta_m-\theta)^T$ and
  $$H(G(x))=(y^T, g(x))^T.$$
Since $ a_1-a, a_2-a, \ldots, a_m-a $   are linearly independent vectors in $\Bbb R^m,$  the matrix $A$ is nonsingular and we have
  $x=A^{-1}y-A^{-1}\hat{\theta}.$ Then
\begin{eqnarray}
	g(x)&=&g\left( A^{-1}y-A^{-1}\hat{\theta}\right)       \label{238qqq}      \\
	&=&y^T  \left( \underbrace{{A^{-1}}^T  A^{-1}}_{B}\right)y-2\underbrace{ \left( \hat{\theta}^T{A^{-1}}^T A^{-1}
		+ a^T{A^{-1}}   \right)}_{\bar{a} ^T}y\\
	& &+\underbrace{\hat{\theta}^T\left( {A^{-1}}^T A^{-1}\right)\hat{\theta}+ 2a ^TA^{-1}\hat{\theta}+\theta}_{\bar{\theta}}\nonumber\\
&:= &\bar{g}(y).\nonumber
\end{eqnarray}
Now we have
$$H(G(x))=(y^T, \bar{g}(y))^T,$$
where  $\bar{g}(y)=y^TBy-2\bar{a}^Ty+\bar{\theta}$ is a strictly convex quadratic function since 
$B={A^{-1}}^T  A^{-1}\succ0.$
The image $H(G(\Bbb R^n))$ is then
\begin{equation}\label{Bk}
	H(G(\Bbb R^n)) =\{(y^T, \bar{g}(y))^T,   y\in\Bbb R^n \}.
\end{equation}
 We apply Lemma \ref{thm21qq} to have
$$ H(G(\R^n))^{\bullet} =\{(y^T, \bar{g}(y))^T,   y\in\Bbb R^n \}^\bullet=\grap(\bar{g})^\bullet=\epi(\bar{g}).$$
Since the function $\bar{g}(y)$ is strictly convex, its epigraph    $\epi(\bar{g})$
is convex and thus $ H(G(\R^n))^\bullet$ is convex. Moreover, $H$ is a nonsingular linear mapping, by Lemma \ref{B0} we have
$ H(G(\R^n))^\bullet= H(G(\R^n)^\bullet).$ This shows that  $H(G(\R^n)^\bullet)$ is convex and 
$ G(\R^n)^\bullet$ is convex.

(ii)  As seen, if  ${\rm rank}\{a_1-a, a_2-a, \cdots, a_m-a\}<n$ then  $G(\R^n)=G(\R^n)^\bullet.$
So $G(\R^n)\cap \Lambda=\emptyset$  implies  $G(\R^n)^\bullet\cap \Lambda=\emptyset.$

 If  ${\rm rank}\{a_1-a, a_2-a, \ldots, a_m-a\}=n$  and  $m=n,$  we will show that
$G(\R^n)\cap \Lambda=\emptyset$  implies $G(\R^n)^\bullet\cap \Lambda=\emptyset.$ Let us assume on contrary that
$G(\R^n)^{\bullet}\cap \Lambda\ne\emptyset.$ Then
there is
\begin{equation}\label{}
	\omega=(\omega_0, \omega_1,\omega_2,\ldots,\omega_n)^T\in G(\R^n)^{\bullet}\cap \Lambda
\end{equation}
    such that   $\omega_0<0, \omega_1\le0, \ldots, \omega_n\le0$
    and $\omega=\lambda u+ (1-\lambda)v$ for some $u\in G(\R^n), v\in G(\R^n),$ $0\le \lambda\le 1.$
Let  
$$T_\omega:=\{\omega+(t,0, \ldots, 0)^T,~ t\leq 0\}=\{(\omega_0+t, \omega_1,\omega_2,\ldots,\omega_n)^T,~ t\leq 0\}\subset \Lambda.$$

We are going to show that $G(\R^n) \cap T_\omega \ne \emptyset,$ which then implies a contradiction that $G(\R^n)\cap \Lambda\ne\emptyset.$

Using the same notations as in (i), we first have, by \eqref{Bk}:
\begin{equation}\label{Bkk}
	H(G(\R^n)) =\{(y^T, \bar{g}(y))^T: y\in \R^n\}.
\end{equation}
Let $z=(y^T, z_0)^T\in\R^{n+1}$ and $f(z)=\bar{g}(y)-z_0,$  then
\begin{equation}\label{Bkkb}
	H(G(\R^n)) =\{z\in \R^{n+1}: f(z)=0 \}.
\end{equation}
Observe that  $f(z)$ is a continuously convex function, its lower level set
\begin{equation}
\Omega:=	\{z \in \R^{n+1}: f(z)\leq 0\}     \label{223nnn}
\end{equation}
is thus a convex set in  $ \R^{n+1}$ and $H(G(\R^n)) \subset\Omega.$
Let  $\bar{u}=H(u),\bar{v}= H(v), \bar{\omega}=H(\omega).$
Since $u\in G(\R^n), v\in G(\R^n)$ and $\omega=\lambda u+ (1-\lambda)v,$ i.e., $\omega$ is in the line segment $[u,v]:$  $\omega\in[u,v],$   we have
$\bar{u}, \bar{v}\in H(G(\R^n))$ and
$\bar{\omega}\in[\bar{u}, \bar{v}].$ This shows $\bar{\omega}\in \Omega$ and thus
\begin{eqnarray}\label{224nn}
	f(\bar{\omega})\le 0.
\end{eqnarray}
 
%Then $\bar{\omega}=L(\omega)=(\omega_1+\omega_{n+1},\omega_2+\omega_{n+1},\ldots,\omega_n+\omega_{n+1}, \omega_{n+1})^T.$
Pick a point  $\omega_t=(\omega_0+t, \omega_1, \ldots, \omega_n)^T\in T_\omega, t\leq 0,$ and consider
\begin{eqnarray}
	 \bar{\omega_t}&=&H(\omega_t)=(\omega_1+\omega_0+t, \ldots, \omega_n+\omega_0+t,- \omega_0-t)^T \nonumber\\
	&=& (\bar{y}_t^T, -\omega_0-t)^T=  (\bar{y}^T+t\bar{1}^T, -\omega_0-t)^T\in H(T_\omega), \label{219bbb}
\end{eqnarray}
  where $ \bar{y}_t^T=(\omega_1+\omega_0+t, \ldots, \omega_n+\omega_0+t),  t\leq 0,
  \bar{y}^T=(\omega_1+\omega_0, \ldots, \omega_n+\omega_0)$ and $ \bar{1}=(1, 1, \ldots, 1)^T,$ we observe that
\begin{eqnarray*} 
	g(t):=f(\bar{\omega}_t)=(\bar{1}^TB\bar{1})t^2+2\left(\bar{1}^TB\bar{y}-\bar{a}^T\bar{1}+\frac{1}{2}\right)t+\omega_0+\bar{g}(\bar{y})
\end{eqnarray*}
is a  quadratic function of $t$ with quadratic coefficient   $\bar{1}^TB\bar{1}>0.$  
So $g(t)>0$ when $t<0$ and $|t|$ is large enough. On the other hand, applying \eqref{224nn} we
have   
$$g(0)=\omega_0+\bar{g}(\bar{y})=f(\bar{\omega})\le0.$$ Those together with the continuity of $g$ imply that there exists $ t^*\leq 0$ such that $g(t^*)=f(\bar{\omega}_{t^*})= 0.$ This indicates that  
  $\bar{\omega}_{t^*}\in  H(G(\R^n)).$ And as seen,  $\bar{\omega}_{t^*}\in H(T_{{\omega}})$ so
   $\bar{\omega}_{t^*}\in H(T_{{\omega}})\cap H(G(\R^n)).$  Since
$H$ is a linear isomorphism, it holds that 
$$H(T_{{\omega}})\cap H(G(\R^n))=H(T_{{\omega}}\cap  G(\R^n)).$$
That is $\bar{\omega}_{t^*}\in H(T_{{\omega}}\cap  G(\R^n))$ and it implies that $ T_{{\omega}} \cap G(\R^n) \ne\emptyset.$ This further implies $\Lambda \cap G(\R^n)\ne \emptyset$ since $T_{{\omega}}\subset \Lambda.$ This contradicts to 
our assumption that $\Lambda \cap G(\R^n)= \emptyset.$
\endproof

%The following example shows that if
%${\rm rank}\{a_1+a_0, a_2+a_0, \cdots, a_m+a_0\}=n$ and $m>n,$ then $G(\R^n)\cap\Lambda=\emptyset$ does not imply $G(\R^n)^{\bullet}\cap\Lambda=\emptyset.$
%\begin{exam}\rm
%Let $g_0(x)=-x^2+2x,$ $g_1(x)=x^2-3,$ $g_2(x)=x^2-4x+3,$ $G=(g_1, g_2, g_0)^T$. Since  $\{g_1(x)\le  0, g_2(x)\le 0, g_0(x)< 0\}$ has no solution, we deduce that  $G(\R^n)\cap\Lambda=\emptyset$. Now, choose $u=0,  v=2$ then $G(u)=(-3,3,0)^T$, $G(v)=(1,-1,0)^T$ and $\frac{1}{4}G(u)+\frac{3}{4}G(v)=(0,0,0)^T\in \Lambda$, it means that $G(\R^n)^{\bullet}\cap\Lambda \ne \emptyset$.
%\end{exam}

As a result we obtain the following extension of the S-Lemma
\cite{Yakubovich,Polik-Terlaky07}.  For completeness, we provide its proof below (c.f. \cite[Theorem 2.2]{Polik-Terlaky07}). 
\begin{cor}[An extension of the S-Lemma]\label{thm2}
	Suppose the interior of the intersection $\cap_{i=1}^mB(a_i,r_i)$ is nonempty, i.e., there is
	a point $\bar{x}\in\R^n$ such that $g_i(\bar{x})<0, i=1,2,\ldots, m.$ Then if either
	${\rm rank}\{a_1-a, a_2-a, \cdots, a_m-a\}<n$ or ${\rm rank}\{a_1-a, a_2-a, \cdots, a_m-a\}=n$ and
	$m=n,$   the following two statements are equivalent.
	\begin{enumerate}
		\item[(i)] The system $\begin{cases} g(x)>0,\\ g_i(x)\le0, i=1,2,\ldots, m, \end{cases}$ is unsolvable,\\
		\item[(ii)] There exist nonnegative scalars $\mu_i\ge0, i=1,2,\ldots,m$ such that
		\begin{eqnarray}\label{B3b}
			 -g(x)+\sum_{i=1}^n\mu_ig_i(x)\ge0 ~\forall x\in\R^n.
		\end{eqnarray}
	\end{enumerate}
\end{cor}
\proof
The implication $(ii)\Rightarrow (i)$ is easy to obtain. We therefore prove only the direction
$(i)\Rightarrow (ii).$ 
%Let  $G(\R^n)$ and
%$\Lambda=\{(z_1, \ldots, z_m, z_{m+1}):  z_i\le0, i=1,2,\ldots,m, z_{m+1}<0\}$ be two sets in $\R^{m+1}$ as introduced in the previous section.
The system (i) is unsolvable meaning that
\begin{align}\label{b}
G(\R^n)\cap\Lambda=\emptyset.
\end{align}
On the other hand, if either
	${\rm rank}\{a_1-a, a_2-a, \cdots, a_m-a\}<n$ or ${\rm rank}\{a_1-a, a_2-a, \cdots, a_m-a\}=n$ and
	$m=n,$ then, by   Theorem \ref{lm222},  $G(\R^n)^\bullet$ is convex and \eqref{b} implies 
\begin{align}\label{bb}
G(\R^n)^\bullet\cap\Lambda=\emptyset.
\end{align}
So, $G(\R^n)^\bullet$ and $\Lambda$ are   separable  by a hyperplane in $\R^{m+1}.$
 That is, there exist scalars $\lambda_0, \lambda_1, \ldots, \lambda_m$ not all zero such that
\begin{align}
	& \lambda_0z_0 +\lambda_1z_1+\ldots+\lambda_mz_m\le 0, \forall z=(z_0, z_1, \ldots, z_m)^T\in\Lambda,\label{H1}\\
	&-\lambda_0g(x)+\lambda_1g_1(x)+\ldots+\lambda_mg_m(x)\ge0,~\forall x\in\R^n.\label{H2}
\end{align}
Since $ (-1,0, \ldots, 0)^T\in\Lambda,$ \eqref{H1} implies that $\lambda_0\ge0.$ Moreover,
substituting
$$u^i=\underbrace{(-\epsilon,0,\ldots, -1,}_{-1 \text{ in the } (i+1)th \text{ position }} 0,\ldots,0)\in\Lambda, i=1,2,\ldots, m,$$
with arbitrary small 
$\epsilon>0$ into \eqref{H1} we obtain that $\lambda_i\ge0$ for $i=1,2,\ldots,m.$
Now, using \eqref{H2} with the assumption that   $g_i(\bar{x})<0$ for all $i=1,2,\ldots,m,$ we see that  $\lambda_0$
cannot be zero. That is  $\lambda_0>0.$
Divide both sides of \eqref{H2} by $\lambda_0$ and let $\mu_i=\frac{\lambda_i}{\lambda_0}$ we get the desired proof.
\endproof

\section{Applications to   the smallest enclosing ball problem}\label{scet3}
Given a set of $m$ balls $B(a_i, r_i),  i=1,2,\ldots,m,$ in $\R^n$ such that the intersection $\cap_{i=1}^mB(a_i,r_i)$ has
a nonempty interior, where  $a_i, r_i $ are the center  and radius
of the ball $B(a_i, r_i).$ We now apply the results in the previous sections
to find a ball $B(a,r)$ of center $a$ and radius $r$ such that $B(a,r)$ is the ball of smallest radius containing
the intersection $\cap_{i=1}^mB(a_i,r_i).$
Note that
$$B(a_i,r_i)=\{x\in\R^n: \|x-a_i\|\le r_i\}, i=1,2,\ldots,m.$$
The enclosing $\cap_{i=1}^mB(a_i,r_i)\subset B(a,r)$ means that
\begin{align}\label{implication}
	\|x-a_i\|\le r_i, i=1,2,\ldots,m, \Rightarrow \|x-a\|\le r.
\end{align}
Equivalently, the system
\begin{align}\label{implication1}
	\begin{cases}
		\|x-a\|^2> r^2,\\
		\|x-a_i\|^2\le r_i^2, i=1,2,\ldots,m,
	\end{cases}
\end{align}
is unsolvable. By Corollary \ref{thm2},
\begin{align}
 \text{\em  if  ~  either ~ } &
{\rm rank}\{a_1-a, a_2-a, \ldots, a_m-a\}<n \text{\em ~  or   }\nonumber\\
&  {\rm rank}\{a_1-a, a_2-a, \ldots, a_m-a\}=n  \text{\em ~ and ~ }
m=n,    \label{note}
\end{align}
 then the un-solvability of \eqref{implication1} is equivalent to the existence of nonnegative real numbers
$\mu_1, \mu_2, \ldots, \mu_m$ such that
\begin{align}\label{SDP1}
	-\|x-a\|^2+r^2+\sum_{i=1}^m\mu_i(\|x-a_i\|^2-r_i^2)\ge 0, ~\forall x\in\R^n.
\end{align}
The inequality \eqref{SDP1} holds if and only if   
\begin{align}\label{SDP2}
	\left(
	\begin{matrix}
		(\sum_{i=1}^m\mu_i-1)I& a-\sum_{i=1}^m\mu_ia_i\\
		a^T-\sum_{i=1}^m\mu_ia_i^T&r^2- \|a\|^2+\sum_{i=1}^m\mu_i(\|a_i\|^2-r_i^2)
	\end{matrix}\right)\succeq0.
\end{align}
Thus, under the assumption \eqref{note},   problem \eqref{problem} is now transformed to the following SDP   problem
of variables $\eta:=r^2>0, a\in\R^n$ and $\mu_i\ge0:$
\begin{align}\label{SDP3}
	{\rm (P)}\hspace*{1.2cm}
	\begin{array}{llll}
		&\min & \eta& \\
		&{\rm s.t.}& \left(
		\begin{matrix}
			(\sum_{i=1}^m\mu_i-1)I& a-\sum_{i=1}^m\mu_ia_i\\
			a^T-\sum_{i=1}^m\mu_ia_i^T&\eta- \|a\|^2+\sum_{i=1}^m\mu_i(\|a_i\|^2-r_i^2)
		\end{matrix}\right)\succeq0,\\
		& & \mu_i\ge0, i=1,2,\ldots,m, a\in\R^n.
	\end{array}
\end{align}

However, $a$ has not been known, we cannot compute the rank of $\{a_1-a, a_2-a, \ldots, a_m-a\}$ to check whether the assumption
 \eqref{note} is satisfied. Fortunately, we always have
${\rm rank}\{a_1-a, a_2-a, \ldots, a_m-a\}\le {\rm rank}\{a_1, a_2, \ldots, a_m\}.$ So we can use 
  ${\rm rank}\{a_1, a_2, \ldots, a_m\}$ instead of ${\rm rank}\{a_1-a, a_2-a, \ldots, a_m-a\}.$ On the other hand,
   the proof of Theorem 3.2 in \cite{Beck07} showed that if $(\eta^*, a^*, \mu^*)$ is an optimal solution 
    of \eqref{SDP3} then $\mu^*=(\mu^*_1, \mu^*_2, \ldots,\mu^*_m)^T$ must satisfy  
     $\sum_{i=1}^m\mu_i^*=1.$  
The matrix inequality in \eqref{SDP3} thus becomes
$$\left(
		\begin{matrix}
			0& a^*-\sum_{i=1}^m\mu^*_ia_i\\
			{a^*}^T-\sum_{i=1}^m\mu^*_ia_i^T&\eta^*- \|a^*\|^2+\sum_{i=1}^m\mu^*_i(\|a_i\|^2-r_i^2)
		\end{matrix}\right)\succeq0.$$
This inequality together with the fact that $\eta^*$ is the smallest number satisfying the inequality indicate that
  $a^*=\sum_{i=1}^m\mu^*_ia_i$ and
$\eta^*=  \|\sum_{i=1}^m \mu_i^* a_i\|^{2}-\sum_{i=1}^m \mu_i^*\left(\|a_i\|^{2}-r_i^{2}\right).$
 As a result, we have the following formula to find the smallest ball $B(a,r).$ The result looks like Theorem 3.2 in \cite{Beck07} but
 under different conditions.

\begin{thm}\label{thm3}
	Let  $B(a_i,r_i),$  $ i=1,2,\ldots,m,$ be the  balls such that the intersection
	$\cap_{i=1}^m B(a_i,r_i)$  has a nonempty interior. If either
	${\rm rank}\{a_1, a_2, \ldots, a_m\}<n$ or   ${\rm rank}\{a_1, a_2, \ldots, a_m\}=n$
	and $n=m,$ then the center $a$ and radius $r$ of the smallest ball $B(a,r)$ enclosing the intersection
$\cap_{i=1}^m B(a_i,r_i)$  are given by
	\begin{eqnarray}\label{}
		a&=&\sum_{i=1}^m \mu_i a_i, \\
		r&=&\sqrt{\|\sum_{i=1}^m \mu_i a_i\|^{2}-\sum_{i=1}^m \mu_i\left(\|a_i\|^{2}-r_i^{2}\right)},
	\end{eqnarray}
	respectively, where $\mu=(\mu_1,\mu_2,\ldots,\mu_m) $ is an optimal solution of the convex quadratic minimization problem
	\[ \begin{array}{l}
		\min \left\|\sum_{i=1}^m \mu_i {a}_i\right\|^{2}-\sum_{i=1}^m \mu_i\left(\left\|{a}_i\right\|^{2}-r_i^{2}\right). \\
		\text { s.t. } \quad \sum_{i=1}^m\mu_i=1, \mu_i\ge0, i=1,2,\ldots,m.
	\end{array}
	\]
	
\end{thm}
 
\section{Conclusion and remarks}

  We have shown that the joint numerical range $G(\R^n)$ of the quadratic mapping   $G=(-g, g_1, \ldots, g_m)$ with $g(x)=x^Tx-2a^Tx+\theta, g_i(x)=x^Tx-2a_i^Tx+\theta_i, i=1,2,\ldots,m,$
is convex if and only if  the rank of $m$ vectors $a_1-a, a_2-a, \ldots, a_m-a\in\Bbb R^n$ is less than $n:$ 
 ${\rm rank}\{a_1-a, a_2-a, \ldots, a_m-a\}<n.$ 
 In case $G(\R^n)$ is not convex, i.e.,  ${\rm rank}\{a_1-a, a_2-a, \ldots, a_m-a\}=n,$ we create an extension $G(\R^n)^\bullet$ of  $G(\R^n)$ in the sense that
   $G(\R^n)\subset G(\R^n)^\bullet,$ and  $G(\R^n)=G(\R^n)^\bullet$ if and only if $G(\R^n)$ is convex. 
   We now obtain a new result that
    if  ${\rm rank}\{a_1-a, a_2-a, \ldots, a_m-a\}=n=m$ then  $G(\R^n)^\bullet$ is convex, even $G(\R^n)$ is not.  
    Interestingly, in this case $G(\R^n)\subsetneq G(\R^n)^\bullet$ but 
    $G(\R^n)\cap \Lambda=\emptyset$ implies  $G(\R^n)^\bullet\cap\Lambda =\emptyset.$
  Those results allow us to obtain a new separable property of 
 $G(\R^n)$   that if $G(\R^n)\cap \Lambda=\emptyset$ then $G(\R^n)$ and $\Lambda$ are separable by a hyperplane if 
 either ${\rm rank}\{a_1-a, a_2-a, \ldots, a_m-a\}<n$ or  ${\rm rank}\{a_1-a, a_2-a, \ldots, a_m-a\}=n=m.$
 As an interesting application, the smallest 
enclosing ball (SEB) problem  of finding a smallest ball $B(a,r)$ in $\R^n$ containing the 
intersection $\cap_{i=1}^mB(a_i, r_i)$ of given balls $B(a_i,r_i)$
 is now solved if 
  either ${\rm rank}\{a_1, a_2, \ldots, a_m\}<n$ or  ${\rm rank}\{a_1, a_2, \ldots, a_m\}=n=m.$ 
  This new progress on solving  SEB problem together with the NP-hard property of  SEB problem \cite{Xia-Yang-Wang20}
   raise a conjecture
that whether the SEB problem  is NP-hard if and only if ${\rm rank}\{a_1, a_2, \ldots, a_m\}=n$ and $m>n?$

\vskip0.3cm
{\bf Declaration of competing interest}

The authors declare that they have no competing interests.
\vskip0.3cm
{\bf Data availability}

No data was used for the research described in the article.

\section*{Funding}

Huu-Quang Nguyen was funded by Vietnam Ministry of Education and Training under Grant number: B2023-TDV-03,
Van-Bong Nguyen
was funded by Tay Nguyen University under Grant number T2025-42CBTD.

%\section*{References}

%\bibliography{mybibfile}

\end{document}